\documentclass[a4paper,12pt,bibliography=totoc,oneside,appendixprefix,fleqn]{scrartcl}
\usepackage{supertabular} 
\usepackage{bbm} 
\usepackage{amsmath} 
\usepackage{amsthm} 
\usepackage{amssymb} 
\usepackage{amscd}
\usepackage{amsfonts}
\usepackage{cancel}
\usepackage{mathcomp}
\usepackage{exscale}
\usepackage{relsize} 
\usepackage{parskip} 
\usepackage{paralist} 
\usepackage[noindentafter]{titlesec} 
\usepackage{booktabs} 
\usepackage{mathtools} 
\usepackage{subfigure} 
\usepackage{siunitx} 
\usepackage{pdfpages} 
\usepackage{alltt} 
\usepackage[margin=10pt,font=small,labelfont=bf,
labelsep=endash]{caption}

\usepackage{latexsym}
\usepackage{latexsym}
\usepackage{enumerate}
\usepackage{verbatim}
\usepackage{color}
\usepackage{graphicx}
\usepackage{tabularx}
\usepackage{arydshln}
\usepackage{booktabs}
\usepackage{euscript}
\usepackage{multicol}
\usepackage{float}

\newtheorem{thm}{Theorem}[section]
\newtheorem{cor}[thm]{Corollary}
\newtheorem{lem}[thm]{Lemma}

\theoremstyle{definition}
\newtheorem{defn}[thm]{Definition}
\newtheorem{rem}[thm]{Remark}
\newtheorem{exa}[thm]{Example}

\newcommand{\N}{\ensuremath{\mathbb{N}}}

\newcommand{\R}{\ensuremath{\mathbb{R}}}

\newcommand{\1}{\ensuremath{\mathbbm{1}}}

\title{Estimating the gap of finite metric spaces of strict p-negative type}
\author{Reinhard Wolf}
\date{\today}

\begin{document}

\setlength{\columnsep}{.01cm}
\begin{multicols}{2}[\section*{Estimating the gap of finite metric spaces of strict p-negative type}
\subsection*{Reinhard Wolf}
\small{Universit\"at Salzburg, Fachbereich Mathematik, Austria}\\
\hrule][.7cm]
ARTICLEINFO\\

\hrule

\small{\emph{Article history:}\\
Received TT month JJJJ\\
Accepted TT month JJJJ\\
Available online TT month JJJJ\\

Submitted by N.N.\\

\hrule

\emph{Keywords:}\\
Finite metric spaces\\
$p$-negative type\\
Matrices of negative type

\begin{abstract}ABSTRACT\\

Let $(X,d)$ be a finite metric space. This paper first discusses the spectrum of the $p$-distance matrix of a finite metric space of $p$-negative type and then gives upper and lower bounds for the so called gap $\Gamma(X,p)$ of a finite metric space of strict $p$-negative type. Furthermore estimations for $\Gamma(X,p)$ under a certain glueing construction for finite metric spaces are given and finally be applied to finite ultrametric spaces.
\end{abstract}}
\end{multicols}

\hrule

\section{Introduction}\label{Intro}

Let $(X,d)$ be a metric space and $p > 0$. Recall that $(X,d)$ has $p$-negative type if for all natural numbers $n$, all $x_1, x_2, \ldots, x_n$ in $X$ and all real numbers $\alpha_1, \alpha_2, \ldots, \alpha_n$ with $\alpha_1 + \alpha_2 + \ldots + \alpha_n = 0$ the inequality
\begin{equation*}
\sum_{i,j = 1}^n \alpha_i \alpha_j d(x_i,x_j)^p \le 0
\end{equation*}
holds.\\
Moreover if $(X,d)$ has $p$-negative type and
\begin{equation*}
\sum_{i,j = 1}^n \alpha_i \alpha_j d(x_i,x_j)^p = 0,\ \text{together with}\ x_i \neq x_j,\ \text{for all}\ i \neq j
\end{equation*}
implies $\alpha_1 = \alpha_2 = \cdots = \alpha_n = 0$, then $(X,d)$ has strict $p$-negative type.

Following \cite{ID1,ID2} we define the $p$-negative type gap $\Gamma(X,p)$ of a $p$-negative type metric space $(X,d)$ as the largest nonnegative constant, such that
\begin{equation*}
\frac{\Gamma(X,p)}{2} \left(\sum_{i = 1}^n |\alpha_i|\right)^2 + \sum_{i,j = 1}^n \alpha_i \alpha_j d(x_i,x_j)^p \le 0
\end{equation*}
holds for all natural numbers $n$, all $x_1, x_2, \ldots, x_n$ in $X$ with $x_i \neq x_j$, for all $i \neq j$, and all real numbers $\alpha_1, \alpha_2, \ldots, \alpha_n$ with $\alpha_1 + \alpha_2 + \ldots + \alpha_n = 0$.

The above defined $p$-negative type gap $\Gamma(X,p)$ can be used to enlarge the $p$-parameter, for wich a given finite metric space is of $p$-negative type: It is shown in \cite{HL} (Theorem 3.3) that a finite metric space $(X,d)$ with cardinality $n = |X| \ge 3$ of strict $p$-negative type is of strict $q$-negative type for all $q \in [p,p + \xi)$ where
\begin{equation*}
\xi = \frac{ln \left(1 + \frac{\Gamma(X,p)}{\Delta(X)^p \cdot \gamma(n)}\right)}{ln \ \EuScript{D}(X)},
\end{equation*}
with 
\begin{equation*}
\begin{multlined}
\Delta(X) = \underset{x,y \in X}{\max} d(x,y),\\
\EuScript{D}(X) = \Delta(X) \cdot \left(\underset{x,y \in X \atop x \neq y}{\min} d(x,y)\right)^{-1} \text{and}\\
\gamma(n) = 1 - \frac{1}{2} \left(\Big{\lfloor} \frac{n}{2}\Big{\rfloor}^{-1} + \Big{\lceil} \frac{n}{2} \Big{\rceil}^{-1}\right).
\end{multlined}
\end{equation*}
For basic information on $p$-negative type spaces ($1$-negative type spaces are also known as quasihypermetric spaces) see for example \cite{ID1,PH,HL,PN1,PN2,PN3,PN4,SS1,AW}.

After section \ref{Intro} and \ref{Notation} (introduction, notation and basic definitions) of this paper we discuss the spectrum of the $p$-distance matrices of finite metric spaces of $p$-negative type in section \ref{Finite} and then focus our attention to upper and lower bounds (section \ref{4},\ref{5}) for the $p$-negative gap $\Gamma(X,p)$ of finite metric spaces of strict $p$-negative type.

In section \ref{6} we estimate $\Gamma(X,p)$ for a certain glueing construction for finite metric spaces and finally apply these results to finite ultrametric spaces in section \ref{7}.

\section{Notation and basic definitions}\label{Notation}

Elements $x$ in $\R^n$ are interpretated as column vectors, so $x^T = (x_1, x_2, \ldots, x_n)$. The canonical inner product of two elements $x, y$ in $\R^n$ is given by $(x|y)$ and the canonical unit vectors are denoted by $e_1, e_2, \ldots, e_n$. The element $\mathbbm{1}$ in $\R^n$ is defined as $\mathbbm{1}^T = (1, 1, \ldots, 1)$.
The usual $p$-norm $(p \ge 1)$ on $\R^n$ is given by $\|x\|_p = \left(\sum_{i = 1}^n |x_i|^p\right)^{\frac{1}{p}}$, for $x^T = (x_1, x_2, \ldots, x_n)$.\\
The subspace $F_0$ and the affine subspace $F_1$ are $F_0 = \{x \in \R^n | (x|\mathbbm{1}) = 0\}$ and\\
$F_1 = \{x \in \R^n | (x|\mathbbm{1}) = 1\}$.

Let $E$ be a subspace of $\R^n$ and $A$ be a real symmetric $n \times n$ matrix. We say $A$ is of nonnegative type (positive type) on $E$ if $(Ax|x) \ge 0$, for all $x$ in $E$ ($(Ax|x) > 0$, for all $x \neq 0$ in $E$).
Hence $A$ is positive semidefinite (positive definite) if and only if $A$ is of nonnegative type (positive type) on $E = \R^n$. We say $A$ is of negative type (strict negative type) on $E$ if $(Ax|x) \le 0$, for all $x$ in $E$ ($(Ax|x) < 0$, for all $x \neq 0$ in $E$).

For a given finite metric space $(X,d)$, $X = \{x_1, x_2, \ldots, x_n\}$ we make the following definitions:\\
The $n \times n$ $p$-distance matrix ($p > 0$) $D_p(X)$ (for short $D_p$, if the underlying space $X$ is clear from context) is given by
\begin{equation*}
D_p(X) (= D_p) = \left(d(x_i,x_j)^p\right)_{i,j = 1}^n.
\end{equation*}

The diameter $\Delta(X)$ (for short $\Delta$, if the underlying space $X$ is clear from context) is defined as
\begin{equation*}
\Delta(X) (= \Delta) = \underset{1 \le i,j \le n}{\max} d(x_i,x_j).
\end{equation*}

For a given $\alpha > 0$ the finite metric space $(X,\alpha d)$ is written as $\alpha X((\alpha d)(x_i,x_j) = \alpha d(x_i,x_j))$. Finally the discrete metric space consisting of $n$ points is denoted by $X_n$.

\section{Finite metric spaces of $p$-negative type and the spectrum of their $p$-distance matrices}\label{Finite}

The following result by Crouzeix and Ferland was used in \cite{TLH1} (Theorem 2.2) to characterize Euclidean distance matrices:\\

\thm \label{Thm1}{(Theorem 4.1 in \cite{IPC}) Let $A$ be a real symmetric $n \times n$ matrix and $a$ in $\R^n$, $a \neq 0$. Then the following are equivalent:
\begin{enumerate}
\item $A$ is of nonnegative type on $E = \{x \in \R^n | (x|a) = 0\}$
\item $A$ is positive semidefinite or $A$ has just one negative eigenvalue and there exists $b$ in $\R^n$ such that $Ab = a$, where $(a|b) \le 0$.
\end{enumerate}}

\normalfont
Let $(X,d)$, $X = \{x_1, x_2, \ldots, x_n\}$ be a finite metric space. By definition $(X,d)$ is of $p$-negative type if and only if $D_p$ is of negative type on $F_0$, i.\,e. $(D_px|x) \le 0$ for all $x$ in $\R^n$ with $(x|\mathbbm{1}) = 0$.\\
If $(X,d)$ has at least two points ($x_i \neq x_j)$ then $(-D_p(e_i + e_j) | e_i + e_j) = -2 d(x_i,x_j)^p < 0$ and hence -$D_p$ is not positive semidefinite.

Applying Theorem \ref{Thm1} to $A = -D_p$ and $a = \mathbbm{1}$ we obtain\\

\thm \label{Thm32}{Let $(X,d)$ be a finite metric space of at least two points and let $\lambda_1 \le \lambda_2 \le \ldots \le \lambda_n$ be the eigenvalues of the $p$-distance matrix $D_p$ of $X$. The following are equivalent:
\begin{enumerate}
\item $(X,d)$ is of $p$-negative type
\item $\lambda_1 \le \lambda_2 \le \ldots \le \lambda_{n-1} \le 0 < \lambda_n$ and there exists some $b$ in $\R^n$ with $D_pb = \mathbbm{1}$ and $(b|\mathbbm{1}) \ge 0$.
\end{enumerate}
}
\normalfont
It was shown in \cite{SS1} that a given finite metric space $(X,d)$ of $p$-negative type of at least two points is of strict $p$-negative type if and only if $D_p$ is non singular and $(D_p^{-1}\mathbbm{1}|\mathbbm{1}) \neq 0$. Combining this result with Theorem \ref{Thm32} we get\\

\thm \label{Thm33}{Let $(X,d)$ be a finite metric space of at least two points and let $\lambda_1 \le \lambda_2 \le \ldots \le \lambda_n$ be the eigenvalues of the $p$-distance matrix $D_p$ of $X$. The following are equivalent:
\begin{enumerate}
\item $(X,d)$ is of strict $p$-negative type
\item $\lambda_1 \le \ldots \le \lambda_{n-1} < 0 < \lambda_n$ and $(D_p^{-1}\mathbbm{1}|\mathbbm{1}) > 0$.
\end{enumerate}}\ \\

\rem \label{Rem4}{Let $(X,d)$ be a finite metric space of $p$-negative type of at least two points. By Theorem \ref{Thm32} there exists some $b$ in $\R^n$ with $D_pb = \mathbbm{1}$ and $(b|\mathbbm{1}) \ge 0$. Note that the value $(b|\mathbbm{1})$ is independent of the choice of $b$:\\
Let $D_pb = D_pb' = \mathbbm{1}$ then $(b|\1) = (b|D_pb') = (D_pb|b') = ( b'|\1)$.}\ \\

\defn{Let $(X,d)$, $X = \{x_1, x_2, \ldots, x_n\}$ be a finite metric space with $p$-distance matrix $D_p$.
\begin{equation*}
M_p(X) (= M_p) := \underset{\alpha_1 + \alpha_2 + \ldots + \alpha_n = 1}{\sup} \sum_{i,j = 1}^n \alpha_i \alpha_j d(x_i,x_j)^p
\end{equation*}
In short form
\begin{equation*}
M_p = \underset{x \in F_1}{\sup} (D_px|x)
\end{equation*}
}

\rem{Let $(X,d)$, $X = \{x_1, x_2, \ldots, x_n\}$ be a finite metric space of at least two points $(x_i \neq x_j)$ then
\begin{equation*}
\left(D_p\left(\frac{e_i + e_j}{2}\right) \Big{|} \frac{e_i + e_j}{2}\right) = \frac{d(x_i,x_j)^p}{2} > 0
\end{equation*}
and therefore $M_p \in (0,\infty]$. For $X = \{x\}$ we clearly have $M_p(\{x\}) = 0$.

The constant $M_p$ was studied in \cite{PN1,PN2,PN3,PN4} for general compact metric spaces of $1$-negative type, it also appears in the study of Euclidean distance matrices (see Corollary 3.2 in \cite{TLH2}). Following the ideas of these papers we obtain\\

\thm\label{Thm37}{(Compare to Theorem 3.1 in \cite{PN1}) Let $(X,d)$ be a finite metric space of at least two points. If $M_p < \infty$ then $(X,d)$ is of $p$-negative type.}
\proof{Assume $(X,d)$ is not of $p$-negative type, then there is some $x_0 \in F_0$ such that $(D_px_0|x_0) > 0$. Choose some $w$ in $F_1$ and let $x_n = nx_0 + w$. Note that $x_n \in F_1$ and $(D_px_n|x_n) = n^2(D_px_0|x_0) + 2n(D_px_0|w) + (D_pw|w)$. Hence $\underset{n \in \N}{\sup}(D_px_n|x_n) = \infty$ and therefore $M_p = \infty$.
\begin{flushright}$\qedsymbol$\end{flushright}}

\thm\label{Thm38}{Let $(X,d)$ be a finite metric space of $p$-negative type of at least two points. Further let $b$ in $\R^n$ with $D_pb = \1$ and $(b|\1) \ge 0$ (see Theorem \ref{Thm32} and Remark \ref{Rem4}). The following are equivalent:
\begin{enumerate}
\item $(b|\1) > 0$
\item $M_p < \infty$
\end{enumerate}
Moreover if $(b|\1) > 0$ then $M_p = (b|\1)^{-1}$.}

\proof{Let $(b|\1) > 0$. Let $u = \frac{b}{(b|\1)}$, so $u \in F_1$ and let $x$ be in $F_1$. Since $(X,d)$ is of $p$-negative type we get $(D_p (x-u) | x-u) \le 0$ and hence $(D_p x|x) \le (b|\1)^{-1}$ and $(D_pu|u) = (b|\1)^{-1}$. Therefore $M_p < \infty$ and $M_p = (b|\1)^{-1}$.

Now let $M_p < \infty$ and assume $(b|\1) = 0$. Fix some $w$ in $F_1$ and let $x_n = w + nb$, for $n \in \N$. Note that $x_n \in F_1$ and $(D_px_n|x_n) = (D_pw|w) + 2n(D_pb|w) + n^2(D_pb|b) = (D_pw|w) + 2n$. Hence $\underset{x \in F_1}{\sup}(D_px|x) = \infty$, a contradiction to $M_p < \infty$.
\begin{flushright}$\qedsymbol$\end{flushright}}

\cor\label{Cor39}{Let $(X,d)$ be a finite metric space of at least two points and let $\lambda_1 \le \lambda_2 \le \ldots \le \lambda_n$ be the eigenvalues of the $p$-distance matrix $D_p$ of $X$. The following are equivalent:
\begin{enumerate}
\item $(X,d)$ is of strict $p$-negative type
\item $\lambda_1 \le \lambda_2 \le \ldots \le \lambda_{n-1} < 0 < \lambda_n$ and $M_p < \infty$.
\end{enumerate}
}

\proof{Combine Theorem \ref{Thm33}, Theorem \ref{Thm37} and Theorem \ref{Thm38}.}

\rem\label{Rem310}{Let $(X,d)$ be a finite metric space of strict $p$-negative type. By Theorem \ref{Thm33} we know that $D_p$ is non singular and $(D_p^{-1}\1|\1) > 0$. If $(X,d)$ consists of at least two points define
\begin{equation*}
u_p := \frac{D_p^{-1}\1}{(D_p^{-1}\1|\1)}
\end{equation*}
and recall $M_p = (D_p^{-1}\1|\1)^{-1}$ by Theorem \ref{Thm38}.
Hence $D_p u_p = M_p\1$ with $u_p$ in $F_1$.

If $(X,d)$ consists of a single point $(X = \{x\})$ we have $D_p u_p = M_p\1$ again, for $u_p = \1$ since $M_p(\{x\}) = 0$.

So in any case there is a uniquely defined element $u_p$ in $F_1$ with $D_p u_p = M_p\1$.
}

\section{An upper bound for $\Gamma(X,p)$}\label{4}
We first need the following

\lem \label{Lem41}{Let $nÊ\in \N$, $n \ge 2$ and $B = (b_{ij})_{i,j = 1}^n$ be a real $n \times n$ matrix with $b_{ii} = 0$, for all $1 \le i \le n$. For $n = 2m$ let $M_1 =  \{x \in \{-1,1\}^n | |\{1Ê\le i \le n | x_i = 1\} | = m\}$ and for $n = 2m + 1$ let $M_2 = \{x \in \{-1 - \frac{1}{m},1\}^n || \{1 \le i \le n | x_i = 1\} | = m + 1\}$. Then we have
\begin{enumerate}
\item For $n = 2m$
\begin{equation*}
\frac{1}{|M_1|}\sum_{x \in M_1} (Bx|x) = (-n) \left(\frac{1}{n(n - 1)} \sum_{i \neq j} b_{ij}\right)
\end{equation*}

\item For $n = 2m + 1$
\begin{equation*}
\frac{1}{|M_2|} \sum_{x \in M_2} (Bx|x) = \left(- \frac{(m + 1)n}{m}\right) \left(\frac{1}{n(n - 1)} \sum_{i \neq j} b_{ij}\right)
\end{equation*}
\end{enumerate}
}

\proof{For $n = 2m$ and $i \neq j$ we get
\begin{equation*}
\sum_{x \in M_1} x_ix_j = 1 \cdot 1 \cdot \binom{2m - 2}{m - 2} + 2 \cdot 1 \cdot (-1) \binom{2m - 2}{m - 1} + (-1)(-1) \binom{2m - 2}{m}
\end{equation*}
and hence
\begin{equation*}
\frac{1}{|M_1|} \sum_{x \in M_1} x_ix_j = - \frac{1}{n - 1}
\end{equation*}

For $n = 2m + 1$ and $i \neq j$ we get
\begin{equation*}
\sum_{x \in M_2} x_ix_j = 1\cdot 1\cdot \binom{2m - 1}{m - 1} + 2\cdot 1\cdot \left(-1 - \frac{1}{m}\right) \binom{2m - 1}{m} + \left(-1 - \frac{1}{m}\right)^2 \binom{2m - 1}{m + 1}
\end{equation*}

and hence
\begin{equation*}
\frac{1}{|M_2|} \sum_{x \in M_2}x_ix_j = -\frac{m + 1}{m(n - 1)}
\end{equation*}

In each case $(n = 2m, n = 2m + 1)$ we have
\begin{equation*}
\frac{1}{|M_k|} \sum_{x \in M_k} (Bx|x) = \sum_{i \neq j}b_{ij} \left(\frac{1}{|M_k|} \sum_{x \in M_k} x_ix_j\right),
\end{equation*}
for $k \in \{1, 2\}$ and so we are done.
\begin{flushright}$\qedsymbol$\end{flushright}
}

\thm \label{Thm42}{Let $(X,d)$, $X = \{x_1, x_2, \ldots, x_n\}$ be a finite metric space of $p$-negative type of at least two points. Then we have
\begin{equation*}
\Gamma(X,p) \le \left(\frac{1}{n(n - 1)} \sum_{i \neq j} d(x_i,x_j)^p\right) \Gamma(X_n,p),
\end{equation*}
where $X_n$ denotes the discrete metric space consisting of $n$ points, where
\begin{equation*}
\Gamma(X_n,p) = \frac{1}{2} \left(\Big{\lfloor} \frac{n}{2}\Big{\rfloor}^{-1} + \Big{\lceil} \frac{n}{2}\Big{\rceil}^{-1}\right) =  \left\{
\begin{aligned}
\frac{2}{n}, \quad n \ \text{even}\\
\frac{2}{n - \frac{1}{n}}, \ n \ \text{odd}
\end{aligned}
\right.
\end{equation*}
}

\proof{By definition of $\Gamma(X,p)$ we have
\begin{equation*}
\Gamma(X,p) \le \frac{2(Bx|x)}{\lVert x\rVert_1^2},
\end{equation*}
for all $x \in F_0$, $x \neq 0$ with $B = -D_p$.

Note  that $M_1$ and $M_2$ defined as in Lemma \ref{Lem41} are both subsets of $F_0$. Hence we get
\begin{equation*}
\Gamma(X,p) \le \frac{1}{|M_k|} \sum_{x \in M_k} \frac{2(Bx|x)}{\lVert x\rVert_1^2},
\end{equation*}
for $k = 1 (n = 2m)$ and $k = 2 (n = 2m + 1)$.

Now $\lVert x\rVert_1 = n$, for all $x$ in $M_1 (n = 2m)$ and $\lVert x\rVert_1 = n + 1$, for all $x$ in $M_2 (n = 2m + 1)$.
So we get
\begin{equation*}
\Gamma(X,p) \le  \left\{
\begin{aligned}
\frac{2}{n^2}\frac{1}{|M_1|} \sum_{x \in M_1}(Bx|x), \quad \ n = 2m\\
\frac{2}{(n + 1)^2}\frac{1}{|M_2|} \sum_{x \in M_2} (Bx|x), \ n = 2m + 1
\end{aligned}
\right.
\end{equation*}

Applying Lemma \ref{Lem41} to $B = -D_p$ and since
\begin{equation*}
\Gamma(X_n,p) = \left\{
\begin{aligned}
\frac{2}{n} \ , \quad \quad n = 2m\\
\frac{2}{n - \frac{1}{n}} , \ n = 2m + 1
\end{aligned}
\right.
\end{equation*}
as shown in \cite{AW}, finishes the proof.
\begin{flushright}$\qedsymbol$\end{flushright}
}

\cor{Let $(X,d)$, $X = \{x_1, x_2, \ldots, x_n\}$ be a finite metric space of $p$-negative type of at least two points. Then we have
\begin{equation*}
\Gamma(X,p) \le \Gamma(\Delta(X)\cdot X_n,p) = \left\{
\begin{aligned}
\Delta(X)^p\cdot \frac{2}{n} \ , n \ \text{even}\\
\Delta(X)^p\cdot \frac{2}{n - \frac{1}{n}}, \ n \ \text{odd}
\end{aligned}
\right.
\end{equation*}

where $X_n$ denotes the discrete space consisting of $n$ points and $\Delta(X)$ denotes the dia\-meter of $X$.

Moreover equality holds if and only if $(X,d)$ is isometrically isomorphic to $\Delta(X)\cdot X_n$.}

\proof{By definition of $\Gamma(X,p)$ and since
\begin{equation*}
\Gamma(X_n,p) = \left\{
\begin{aligned}
\frac{2}{n} \ , \quad \quad n \ \text{even}\\
\frac{2}{n - \frac{1}{n}} , \ n \ \text{odd}
\end{aligned}
\right.
\end{equation*}
as shown in \cite{AW} we have

\begin{equation*}
\Gamma(\Delta(X)\cdot X_n,p) = \Delta(X)^p \cdot \Gamma(X_n,p) = \left\{
\begin{aligned}
\Delta(X)^p\cdot \frac{2}{n} \ , \quad \quad n \ \text{even}\\
\Delta(X)^p\cdot \frac{2}{n - \frac{1}{n}} , \ n \ \text{odd}.
\end{aligned}
\right.
\end{equation*}

Now assume that $(X,d)$ is not isometrically isomorphic to $\Delta(X)\cdot X_n$. Hence there exists $i \neq j$ such that $d(x_i,x_j) < \Delta(X)$. By Theorem \ref{Thm42} we get
\begin{equation*}
\Gamma(X,p) < \Delta(X)^p\cdot \Gamma(X_n,p).
\end{equation*}
\begin{flushright}$\qedsymbol$\end{flushright}
}

\section{Bounds for $\Gamma(X,p)$ involving the spectrum}\label{5}

\thm\label{Thm51}{Let $(X,d)$, $X = \{x_1, x_2, \ldots, x_n\}$ be a finite metric space of strict $p$-negative type of at least two points. Let $\lambda_1 \le \lambda_2 \le \ldots \le \lambda_{n-1} < 0 < \lambda_n$ be the eigenvalues of the $p$-distance matrix $D_p$ of $X$ (see Theorem \ref{Thm33}). Then we have
\begin{equation*}
\left(\frac{\lambda_n}{nM_p}\right)|\lambda_{n-1}| \Gamma(X_n,p) \le \Gamma(X,p) \le |\lambda_{n-1}|,
\end{equation*}
where $X_n$ denotes the discrete metric space consisting of $n$ points, where
\begin{equation*}
\Gamma(X_n,p) = \left\{
\begin{aligned}
\frac{2}{n} \ , \ n \ \text{even}\\
\frac{2}{n - \frac{1}{n}} \ , \ n \ \text{odd}
\end{aligned}
\right.
\end{equation*}
and $M_p = (D_p^{-1}\1|\1)^{-1}$ (see Theorem \ref{Thm33} and Theorem \ref{Thm38})
}
\proof{The upper bound:\\
We first show, that $\lVert x\rVert_2^2 \le \frac{1}{2}$ , for all $x$ in $F_0$ with $\lVert x\rVert_1 = 1$. Assume (w.\,l.\,o.\,g.) that $x$ in $F_0$ with  $\lVert x\rVert_1 = 1$ is given by $x = (x_1, \ldots, x_k, -x_{k + 1}, \ldots, -x_n)$, with $x_1, x_2, \ldots, x_k > 0, x_{k+1}, \ldots, x_n \le 0$ and $1 \le k \le n - 1$.

Since $x \in F_0$ and $\lVert x\rVert_1 = 1$ we get $x_1 + \ldots + x_k = x_{k+1} + \ldots + x_n = \frac{1}{2}$. Hence $\lVert x\rVert_2^2 = (x_1^2 + \ldots + x_k^2) + (x_{k+1}^2 + \ldots + x_n^2) \le (x_1 + \ldots + x_k)^2 + (x_{k+1} + \ldots + x_n)^2 = \frac{1}{2}$.

Note that the eigenvalues of $-D_p$ are $-\lambda_n < 0 < |\lambda_{n-1}| \le \ldots \le |\lambda_1|$. Applying the well known Courent-Fisher Theorem we get
\begin{equation*}
|\lambda_{n-1}| = \underset{E}{\sup}\left(\underset{x \in E \atop x \neq 0}{\inf} \frac{(-D_px|x)}{\lVert x\rVert_2^2}\right),
\end{equation*}
where the supremum is taken over all $n - 1$ dimensional subspaces $E$ of $\R^n$. Therefore
\begin{equation*}
|\lambda_{n-1}| \ge \underset{x \in F_0 \atop x \neq 0}{\inf} \frac{(-D_px|x)}{\lVert x\rVert_2^2}.
\end{equation*}

Since $\lVert x\rVert_2^2 \le \frac{1}{2}$, for all $x$ in $F_0$ with $\lVert x\rVert_1 = 1$ (as shown above) we obtain
\begin{equation*}
|\lambda_{n-1}| \ge \underset{x \in F_0 \atop x \neq 0}{\inf} \frac{(-D_px|x)}{\lVert x\rVert_1^2} \cdot \frac{1}{\lVert \frac{x}{\lVert x\rVert_1}\rVert_2^2} \ge 2\underset{x \in F_0 \atop x \neq 0}{\inf} \frac{(-D_px|x)}{\lVert x\rVert_1^2}.
\end{equation*}

Since (by definition)
\begin{equation*}
\frac{\Gamma(X,p)}{2} = \underset{x \in F_0 \atop x \neq 0}{\inf} \frac{(-D_px|x)}{\lVert x\rVert_1^2} ,
\end{equation*}
we get $|\lambda_{n-1}| \ge \Gamma(X,p)$.

The lower bound:\\
Let $\{f_1, f_2, \ldots, f_n\}$ be an orthonormal basis consisting of eigenvectors of $D_p$ ($D_pf_i = \lambda_if_i$, $1 \le i \le n$).\\
By the well known Frobenius - Perron Theorem we know that all coordinates of $f_n$ are nonnegative and hence $(\1|f_n) > 0$. For $x$ in $F_0$ we get $0 = (x|1) = \sum_{i = 1}^n (x|f_i)(\1|f_i)$ and therefore $(x|f_n)^2 = \frac{1}{(\1|f_n)^2} \left(\sum_{i = 1}^{n - 1}(x|f_i)(\1|f_i)\right)^2$

Since

\begin{equation*}
\begin{multlined}
\left(\sum_{i = 1}^{n - 1}(x|f_i)(\1|f_i)\right)^2 = \left(\sum_{i = 1}^{n - 1}|\lambda_i|^{\frac{1}{2}}(x|f_i)|\lambda_i|^{-\frac{1}{2}}(\1|f_i)\right)^2 \\
\le \left(\sum_{i = 1}^{n - 1}|\lambda_i|(x|f_i)^2\right) \cdot \left(\sum_{i = 1}^{n - 1}|\lambda_i|^{-1}(\1|f_i)^2\right)
\end{multlined}
\end{equation*}
 we obtain
 
 \begin{equation*}
 (x|f_n)^2 \le Ê\frac{1}{(\1|f_n)^2}\left(\sum_{i = 1}^{n -1}|\lambda_i|(x|f_i)^2\right)\cdot \left(\sum_{i = 1}^{n - 1}|\lambda_i|^{-1}(\1|f_i)^2\right)
 \end{equation*}
 
 It follows that
 
\begin{equation*}
\begin{multlined}
(D_px|x) = \lambda_n(x|f_n)^2 - \sum_{i = 1}^{n - 1}|\lambda_i|(x|f_i)^2 \\
\le \left(\sum_{i = 1}^{n - 1}|\lambda_i|(x|f_i)^2\right) \cdot \left(\frac{\lambda_n}{(\1|f_n)^2}\sum_{i = 1}^{n - 1}|\lambda_i|^{-1}(\1|f_i)^2 - 1\right)
\end{multlined}
\end{equation*}

Since $(D_p^{-1}\1|\1) = \lambda_n^{-1}(\1|f_n)^2 - \sum_{i = 1}^{n - 1}|\lambda_i|^{-1}(\1|f_i)^2$ we obtain
\begin{equation*}
(-D_px|x) \ge \frac{\lambda_n}{(\1|f_n)^2}(D_p^{-1}\1|\1)\left(\sum_{i = 1}^{n - 1}|\lambda_i|(x|f_i)^2\right)
\end{equation*}
Recall that
\begin{equation*}
\begin{multlined}
(x|f_n)^2 = \frac{1}{(\1|f_n)^2}\left(\sum_{i = 1}^{n - 1}(x|f_i)(\1|f_i)\right)^2\\
\le \frac{1}{(\1|f_n)^2}\left(\sum_{i = 1}^{n - 1}(x|f_i)^2\right)\left(\sum_{i = 1}^{n - 1}(\1|f_i)^2\right).
\end{multlined}
\end{equation*}

Now
\begin{equation*}
\begin{multlined}
\lVert x\rVert_2^2 = (x|f_n)^2 + \sum_{i = 1}^{n - 1}(x|f_i)^2 \le \\
\le \left(\sum_{i = 1}^{n - 1}(x|f_i)^2\right)\left(\frac{1}{(\1|f_n)^2}\sum_{i = 1}^{n - 1}(\1|f_i)^2 + 1\right) = \\
= \frac{n}{(\1|f_n)^2}\sum_{i = 1}^{n - 1}(x|f_i)^2
\end{multlined}
\end{equation*}
and hence
\begin{equation*}
\sum_{i = 1}^{n - 1}(x|f_i)^2 \ge \frac{(\1|f_n)^2}{n}\lVert x\rVert_2^2.
\end{equation*}

Summing up we get
\begin{equation*}
\begin{multlined}
(-D_px|x) \ge \frac{\lambda_n}{(\1|f_n)^2}(D_p^{-1}\1|\1)\left(\sum_{i = 1}^{n - 1}|\lambda_i|(x|f_i)^2\right) \ge \\
\ge \frac{\lambda_n}{(\1|f_n)^2}(D_p^{-1}\1|\1)|\lambda_{n - 1}|\left(\sum_{i = 1}^{n - 1}(x|f_i)^2\right) \ge \\
\ge \frac{\lambda_n}{(\1|f_n)^2}(D_p^{-1}\1|\1)|\lambda_{n - 1}| \cdot \frac{(\1|f_n)^2}{n}\lVert x\lVert_2^2
\end{multlined}
\end{equation*}
So it follows that
\begin{equation*}
(-D_px|x) \ge \frac{|\lambda_{n - 1}|\lambda_n}{nM_p}\lVert x\rVert_2^2,
\end{equation*}
for all $x$ in $F_0$.

By definition of $\Gamma(X,p)$ we have
\begin{equation*}
\Gamma(X,p) = 2\cdot \underset{x \in F_0 \atop x \neq 0}{\inf} \frac{(-D_px|x)}{\lVert x\rVert_1^2} = 2\cdot\underset{x \in F_0 \atop \lVert x\rVert_1 = 1}{\inf}(-D_px|x)
\end{equation*}

Therefore
\begin{equation*}
\Gamma(X,p) \ge \frac{2|\lambda_{n - 1}|\lambda_n}{nM_p}\cdot\underset{x \in F_0 \atop \lVert x\rVert_1 = 1}{\inf} \lVert x\rVert_2^2.
\end{equation*}

By definition of $\Gamma(X_n,p)$ we have
\begin{equation*}
\begin{multlined}
\Gamma(X_n,p) = 2\cdot \underset{x \in F_0 \atop \lVert x\rVert_1 = 1}{\inf}(-D_p(X_n)x|x) \\
= 2\cdot \underset{x \in F_0 \atop \lVert x\rVert_1 = 1}{\inf}\left(\left(\sum_{i = 1}^n e_ie_i^T - \1 \1^T\right)(x)|x\right) = 2\cdot \underset{xÊ\in F_0 \atop \lVert x\rVert_1 = 1}{\inf}\lVert x\rVert_2^2
\end{multlined}
\end{equation*}
and hence we are done.
\begin{flushright}$\qedsymbol$\end{flushright}
}

\rem\label{Rem52}{The factor $\frac{\lambda_n}{nM_p}$ appearing in Theorem \ref{Thm51} is less or equal to one and one if and only if $D_p$ has constant row sum:
$(D_p^{-1}\1|\1)^{-1} = M_p = \underset{x \in F_1}{\sup}(D_px|x)$ and as noted in the proof of Theorem \ref{Thm51} we have $(\1|f_n) > 0$ and hence
\begin{equation*}
M_p \ge \frac{1}{(\1|f_n)^2}(D_pf_n|f_n) = \frac{\lambda_n}{(\1|f_n)^2}
\end{equation*}

Therefore
\begin{equation*}
\frac{\lambda_n}{nM_p} \le \frac{(\1|f_n)^2}{n} \le 1 \quad \text{and} \quad \frac{\lambda_n}{nM_p} = 1
\end{equation*}
if and only if $f_n = \frac{1}{\sqrt{n}}\1$, wich is equivalent to $D_p$ has constant row sum.

\cor{Let $(X,d)$, $X = \{x_1, x_2, \ldots, x_n\}$ be a finite metric space of strict $p$-negative type of at least two points. Let $\lambda_1 \le \lambda_2 \le \ldots \le \lambda_{n - 1} < 0 < \lambda$ be the eigenvalues of the $p$-distance matrix $D_p$ of $X$ and assume that $D_p$ has constant row sum. Then we have
\begin{equation*}
\Gamma(X,p) \ge |\lambda_{n - 1}|\Gamma(X_n,p),
\end{equation*}
where $X_n$ denotes the discrete metric space consisting of $n$ points, where
\begin{equation*}
\Gamma(X_n,p) = \left\{
\begin{aligned}
\frac{2}{n} \ , \ n \quad \text{even}\\
\frac{2}{n - \frac{1}{n}} \ , \ n \quad \text{odd}.
\end{aligned}
\right.
\end{equation*}
 Moreover the inequality is sharp for $X = \Delta(X)X_n$.
 }
 \proof{Combine Theorem \ref{Thm51} and Remark \ref{Rem52} and note that the eigenvalues of the $p$-distance matrix of $\Delta(X)X_n$ are $\lambda_1 = \lambda_2 = \ldots = \lambda_{n - 1} = -\Delta(X)^p$ and\\
 $\lambda_n = (n - 1)\Delta(X)^p$. Hence
 \begin{equation*}
 \Gamma(\Delta(X)X_n,p) = \Delta(X)^p\Gamma(X_n,p) = |\lambda_{n - 1}|\cdot \Gamma(X_n,p).
 \end{equation*}
 \begin{flushright}$\qedsymbol$\end{flushright}
 }

\section{Estimating $\Gamma(X,p)$ for a certain glueing construction}\label{6}

In \cite{SS2} S\'{a}nchez showed (Theorem 4.8) that $\Gamma(X,p)^{-1}$ behaves additively for an additive combination of finite metric spaces of strict $p$-negative type, where an additive combination of two metric spaces is, roughly speaking, a space made by picking a point in each space and glueing them together.

We now consider another glueing construction often used in \cite{PN2,PN3,PN4}.
\defn{(see Theorem 3.5 of \cite{PN2})Let $(X_1,d_1)$ and $(X_2,d_2)$ be two finite metric spaces such that $X_1 \cap X_2 = \emptyset$. Further let $c > 0$ with $2c \ge \max(\Delta(X_1),\Delta(X_2))$. We define
\begin{equation*}
X := X_1 \cup X_2
\end{equation*}
and $d : X \times X \rightarrow \R$ as
\begin{equation*}
d(x,y) := \left\{
\begin{aligned}
d_1(x,y), \ \text{for} \ x, y \in X_1\\
d_2(x,y), \ \text{for} \ x, y \in X_2\\
c, \ \text{for} \ x \in X_1 \ \text{and} \ y \in X_2\\
c, \ \text{for} \ x \in X_2 \ \text{and} \ y \in X_1
\end{aligned}
\right.
\end{equation*}

The space $(X,d)$ will be denoted by $(X_1cX_2,d)$.

The following result was shown in \cite{PN2} (Theorem 3.5) for general compact metric spaces of $1$-negative type. For completeness we add the proof now given in terms of matrices.\\
Note that Corollary \ref{Cor39} implies that $M_p(X_1)$ and $M_p(X_2)$ are both finite, if $(X_1,d_1)$ and $(X_2,d_2)$ are of strict $p$-negative type.

\thm{Let $(X_1,d_1)$ and $(X_2,d_2)$ be two finite metric spaces and $c > 0$ such that $2c \ge \max(\Delta(X_1),\Delta(X_2))$. Then we have
\begin{enumerate}
\item $(X_1cX_2,d)$ is a finite metric space.
\item If $(X_1,d_1)$ and $(X_2,d_2)$ are of $p$-negative type and $M_p(X_1)$ and $M_p(X_2)$ are both finite, then $(X_1cX_2,d)$ is of $p$-negative type if and only if $2c^p \ge M_p(X_1) + M_p(X_2)$.
\item If $(X_1,d_1)$ and $(X_2,d_2)$ are of strict $p$-negative type then $(X_1cX_2,d)$ is of strict $p$-negative type if and only if $2c^p > M_p(X_1) + M_p(X_2)$.
\end{enumerate}
}
\proof{
\begin{description}
\item[\normalfont{assertion 1:}]  clear, since $2c \ge \max(\Delta(X_1),\Delta(X_2))$
\item[\normalfont{assertion 2:}] Let $|X_1| = n$ and $|X_2| = m$. By definition of $X_1cX_2$ the $p$-distance matrix $D_p(X_1cX_2)$ of $X_1cX_2$ is the $(n + m)\times(n + m)$ matrix given by
\begin{equation*}
D_p = D_p(X_1cX_2) = \left[D_p(X_1) \quad c^p\1\1^T \atop c^p\1\1^T \quad D_p(X_2)\right]
\end{equation*}
($\1$ denotes the vector consisting of all ones in arbitraty dimensions) Let $z \in \R^{n + m}$, $z^T = (x^T,y^T)$ with $x \in \R^n$ and $y \in \R^m$. We get
\begin{equation*}
(D_pz|z) = (D_p(X_1)x|x) + (D_p(X_2)y|y) + 2c^p(x|\1)(y|\1).
\end{equation*}

Assume that $(X_1cX_2,d)$ is of $p$-negative type and $2c^p < M_p(X_1) + M_p(X_2)$. By definition of $M_p(X_1)$ and $M_p(X_2)$ there are $x_1 \in F_1$, $x_2 \in F_1$ such that\\
$2c^p < (D_p(X_1)x_1|x_1) + (D_p(X_2)x_2|x_2)$. Let $z = (x_1^T,-x_2^T)^T$. Note that\\
$(z|\1) = (x_1|\1) - (x_2|\1) = 0$.\\
Now $(D_pz|z) = (D_p(X_1)x_1|x_1) + (D_p(X_2)x_2|x_2) - 2c^p(x_1|\1)(x_2|\1)$ and hence $(D_pz|z) = (D_p(X_1)x_1|x_1) + (D_p(X_2)x_2|x_2) - 2c^p > 0$, a contradiction to $(X_1cX_2,d)$ is of $p$-negative type.

On the other hand let $2c^p \ge M_p(X_1) + M_p(X_2)$ and assume $z \in F_0(z^T = (x^T,y^T))$.
Hence $0 = (z|\1) = (x|\1) + (y|\1)$ and assume first that $(x|\1) = 0$. Then we get $(y|\1) = 0$ and $(D_pz|z) = (D_p(X_1)x|x) + (D_p(X_2)y|y) \le 0$, since $(X_1,d_1)$ and $(X_2,d_2)$ are of $p$-negative type.\\
If $(x|\1) \neq 0$ we get $(y|\1) = -(x|\1)$ and
\begin{equation*}
(D_pz|z) = (x|\1)^2\left(D_p(X_1)\frac{x}{(x|\1)} \Big{|} \frac{x}{(x|\1)}\right) + (x|\1)^2\left(D_p(X_2)\frac{y}{(y|\1)} \Big{|} \frac{y}{(y|\1)}\right) - 2c^p(x|\1)^2.
\end{equation*}

Therefore
\begin{equation*}
(D_pz|z) \le (x|\1)^2(M_p(X_1) + M_p(X_2) - 2c^p) \le 0.
\end{equation*}

In each case we have $(D_pz|z) \le 0$, for all $z \in F_0$ and hence $(X_1cX_2,d)$ is of $p$-negative type.
\item[\normalfont{assertion 3:}] Assume that $(X_1cX_2,d)$ is of strict $p$-negative type and $2c^p = M_p(X_1) + M_p(X_2)$. By Remark \ref{Rem310} we get $u_1 \in F_1$ and $u_2 \in F_2$ such that $D_p(X_1)u_1 = M_p(X_1)\1$ and $D_p(X_2)u_2 = M_p(X_2)\1$. Let $z = (u_1^T,-u_2^T)^T$. Note that $z \neq 0$ and $(z|\1) = (u_1|\1) - (u_2|\1) = 0$.
Now
$(D_pz|z) = M_p(X_1) + M_p(X_2) - 2c^p = 0$, a contradiction to $(X_1cX_2,d)$ is of strict $p$-negative type.\\
On the other hand let $2c^p > M_p(X_1) + M_p(X_2)$ and assume $z \in F_0(z^T = (x^T,y^T))$ such that $(D_pz|z) = 0$.\\
If $(x|\1) = 0$ we get $(y|\1) = 0$ and so $0 = (D_pz|z) = (D_p(X_1)x|x) + (D_p(X_2)y|y) \le 0$ implies $x = 0$ and $y = 0$ and hence $z = 0$, since $(X_1,d_1)$ and $(X_2,d_2)$ are of strict $p$-negative type.
If $(x|\1) \neq 0$ we get $(y|\1) = -(x|\1)$ and so
\begin{equation*}
0 = (x|\1)^2\left(D_p(X_1)\frac{x}{(x|\1)} \Big{|} \frac{x}{(x|\1)}\right) + (x|\1)^2\left(D_p(X_2)\frac{y}{(y|\1)} \Big{|} \frac{y}{(y|\1)}\right) - 2c^p(x|\1)^2
\end{equation*}
implies

\begin{equation*}
\begin{multlined}
0 = \left(D_p(X_1)\frac{x}{(x|\1)} \Big{|} \frac{x}{(x|\1)}\right) + \left(D_p(X_2)\frac{y}{(y|\1)} \Big{|} \frac{y}{(y|\1)}\right) - 2c^p \\
\le M_p(X_1) + M_p(X_2) - 2c^p < 0 ,
\end{multlined}
\end{equation*}
a contradiction.
\end{description}

 \begin{flushright}$\qedsymbol$\end{flushright}
}

In order to estimate $\Gamma_p(X_1cX_2)$ we need the following (the proof is straight forward)

\lem\label{Lem63}{
\begin{enumerate}
\item Let $A_1, A_2$ be two real symmetric invertible $n \times n$, $m \times m$ matrices and $c$ be a real number with $1 \neq c^2(x_1|\1)(x_2|\1)$, where $x_1 = A_1^{-1}\1$ and $x_2 = A_2^{-1}\1$.
The inverse matrix of the $(n + m) \times (n + m)$ matrix

$B = \begin{bmatrix} A_1 & c\1\1^T\\ c\1\1^T & A_2\end{bmatrix} $is given by
\begin{equation*}
B^{-1} = \begin{bmatrix} A_1^{-1} + \alpha x_1x_1^T & \beta x_1x_2^T \\ \beta x_2x_1^T & A_2^{-1} + \gamma x_2x_2^T \end{bmatrix},
\end{equation*}
where
\begin{equation*}
\begin{multlined}
\alpha = c^2(x_2|\1)(1 - c^2(x_1|\1)(x_2|\1))^{-1},\\
\beta = -c(1 - c^2(x_1|\1)(x_2|\1))^{-1} \ \text{and}\\
\gamma = c^2(x_1|\1)(1 - c^2(x_1|\1)(x_2|\1))^{-1}.
\end{multlined}
\end{equation*}
\item Let $A$ be a real symmetric invertible $n \times n$ matrix with $(x|\1) \neq 0$, where $x = A^{-1}\1$. Further let $c$ be a real number, $c \neq 0$.
The inverse matrix of the $(n + 1)\times(n + 1)$ matrix
$B = \begin{bmatrix} A & c\1 \\ c\1^T & 0\end{bmatrix}$ is given by
\begin{equation*}
B^{-1} = \begin{bmatrix} A^{-1} - \frac{xx^T}{(x|\1)} & \frac{1}{c(x|\1)} x \\ \frac{1}{c(x|\1)} x^T & -\frac{1}{c^2(x|\1)} \end{bmatrix}
\end{equation*}
\end{enumerate}
}\ \\

\cor\label{Cor64}{\begin{enumerate}
\item Let $(X_1,d_1)$, $(X_2,d_2)$; $|X_1| \ge 2$, $|X_2| \ge 2$ be two finite me\-tric spaces of strict $p$-negative type. Further let $c$ be a real number with $2c \ge \max(\Delta(X_1),\Delta(X_2))$ and $2c^p > M_p(X_1) + M_p(X_2)$. Then we have:\\
The inverse matrix of $D_p(X_1cX_2)$ is given by
\begin{equation*}
D_p(X_1cX_2)^{-1} = \begin{bmatrix} D_p(X_1)^{-1} + \alpha x_1x_1^T & \beta x_1x_2^T \\ \beta x_2x_1^T & D_p(X_2)^{-1} + \gamma x_2x_2^T\end{bmatrix},
\end{equation*}
where
\begin{equation*}
\begin{multlined}
x_1 = D_p(X_1)^{-1}\1, x_2 = D_p(X_2)^{-1}\1,\\
\alpha = c^{2p}(x_2|\1)(1 - c^{2p}(x_1|\1)(x_2|\1))^{-1},\\
\beta = -c^p(1 - c^{2p}(x_1|\1)(x_2|\1))^{-1} \text{and}\\
\gamma = c^{2p}(x_1|\1)(1 - c^{2p}(x_1|\1)(x_2|\1))^{-1}.
\end{multlined}
\end{equation*}

\item Let $(X_1,d_1)$, $(X_2,d_2)$, $|X_1| \ge 2$, $|X_2| = 1$ be two finite metric spaces of strict $p$-negative type. Further let $c$ be a real number with $2c \ge \Delta(X_1)$ and $2c^p > M_p(X_1)$ (Recall $M_p(X_2) = 0$).

The inverse matrix of $D_p(X_1cX_2)$ is given by
\begin{equation*}
D_p(X_1cX_2)^{-1} = \begin{bmatrix} D_p(X_1)^{-1} - \frac{x_1x_1^T}{(x_1|\1)} & \frac{1}{c^p(x_1|\1)}x_1 \\ \frac{1}{c^p(x_1|\1)}x_1^T & -\frac{1}{c^{2p}(x_1|\1)}\end{bmatrix},
\end{equation*}
where $x_1 = D_p(X_1)^{-1}\1$.
\end{enumerate}
}

\proof{\begin{description}
\item[\normalfont{assertion 1:}] Since $M_p(X_1)M_p(X_2) \le \left(\frac{M_p(X_1) + M_p(X_2)}{2}\right)^2$ we get $M_p(X_1)M_p(X_2) < c^{2p}$.
Recall that $M_p(X_1) = (x_1|\1)^{-1}$ and $M_p(X_2) = (x_2|\1)^{-1}$ and  hence we are done by Lemma \ref{Lem63} part 1.
\item[\normalfont{assertion 2:}] Since $(D_p(X_1)^{-1}\1|\1) > 0$, we can apply Lemma \ref{Lem63} part 2.
\end{description}
\begin{flushright}$\qedsymbol$\end{flushright}
}

Let $A$ be a real symmetric invertible $n \times n$ matrix with $(A^{-1}\1|\1) \neq 0$. We define the $n \times n$ matrix $\hat{A}$ by
\begin{equation*}
\hat{A} = \frac{(A^{-1}\1)(A^{-1}\1)^T}{(A^{-1}\1|\1)} - A^{-1}.
\end{equation*}

\lem\label{Lem65}{(The proof is longer, but straight forward)

\begin{enumerate}
\item Let $A_1, A_2, c, x_1, x_2$ and $B$ be defined as in Lemma \ref{Lem63} part 1. Further assume that $(x_1|\1) \neq 0$, $(x_2|\1) \neq 0$ and $2c \neq (x_1|\1)^{-1} + (x_2|\1)^{-1}$.

Then we have
\begin{equation*}
\begin{multlined}
(\hat{B}z|z) = (\hat{A}_1x|x) + (\hat{A}_2y|y) + \\
+ (2c - (x_1|\1)^{-1} - (x_2|\1)^{-1})^{-1}\cdot\left(\left(\frac{x_1}{(x_1|\1)} \Big{|} x\right) - \left(\frac{x_2}{(x_2|\1)} \Big{|} y\right)\right)^2,
\end{multlined}
\end{equation*}
for all $z = (x^T,y^T)^T$, $x$ in $\R^n$ and $y$ in $\R^m$.

\item Let $A, c, x$ and $B$ be defined as in Lemma \ref{Lem63} part 2. Further assume that $2c \neq (x|\1)^{-1}$. Then we have
\begin{equation*}
(\hat{B}z|z) = (\hat{A}u|u) + (2c - (x|\1	)^{-1})^{-1}\left(\left(\frac{x}{(x|\1)} \Big{|} u\right) - v\right)^2,
\end{equation*}
for all $z = (u^T,v)^T$, $u$ in $\R^n$ and $v$ in $\R$.
\end{enumerate}
}\ \\

\cor\label{Cor66}{\begin{enumerate}
\item Let $(X_1,d_1)$, $(X_2,d_2)$, $|X_1| \ge 2$, $|X_2| \ge 2$ be two finite me\-tric spaces of strict $p$-negative type. Further let $c$ be a real number with\\
$2c \ge \max(\Delta(X_1),\Delta(X_2)$ and $2c^p > M_p(X_1) + M_p(X_2)$. Then we have
\begin{equation*}
\begin{multlined}
\widehat{(D_p(X_1cX_2)}z|z) = (\widehat{D_p(X_1)}x|x) + (\widehat{D_p(X_2)}y|y) + \\
+ (2c^p - M_p(X_1) - M_p(X_2))^{-1}\cdot((u_p^1|x) - (u_p^2|y))^2,
\end{multlined}
\end{equation*}
for all $z = (x^T,y^T)^T$, $x$ in $\R^n$, $y$ in $\R^m$, where $D_p(X_1)u_p^1 = M_p(X_1)\1$,\\
$D_p(X_2)u_p^2 = M_p(X_2)\1$,\\
$u_p^1, u_p^2$ in $F_1$ (see Remark \ref{Rem310}).

\item Let $(X_1,d_1)$, $(X_2,d_2)$, $|X_1| \ge 2$, $|X_2| = 1$ be two finite metric spaces of strict

$p$-negative type. Further let $c$ be a real number with $2c \ge \Delta(X_1)$ and $2c^p > M_p(X_1)$. Then we have
\begin{equation*}
\widehat{(D_p(X_1cX_2)}z|z) = \widehat{(D_p(X_1)}u|u) + (2c^p - M_p(X_1))^{-1}((u_p^1|u) - v)^2,
\end{equation*}
for all $z = (u^T,v)^T$, $u$ in $\R^n$, $v$ in $\R$, where $D_p(X_1)u_p^1 = M_p(X_1)\1$, $u_p^1$ in $F_1$ (see Remark \ref{Rem310}).
\end{enumerate}
}

\proof{The assertions follow by Corollary \ref{Cor64} and Lemma \ref{Lem65}.
\begin{flushright}$\qedsymbol$\end{flushright}
}

\thm\label{Thm67}{Let $(X_1,d_1)$, $(X_2,d_2)$, $\max(|X_1|, |X_2|) \ge 2$, be two finite metric spaces of strict $p$-negative type. Let $u_p^1, u_p^2$ be the uniquely defined elements in $F_1$ such that $D_p(X_k)u_p^k = M_p(X_k)\1$, $k = 1,2$, where $M_p(X_k) = \underset{x \in F_1}{\sup}(D_p(X_k)x|x)$ and $M_P(X_k) = (D_p(X_k)^{-1}\1|\1)^{-1}$, if $|X_k| \ge 2$. Further let $c$ be a real number with $2c \ge \max(\Delta(X_1),\Delta(X_2))$ and $2c^p > M_p(X_1) + M_p(X_2)$. Then we have

\begin{equation*}
\left(\frac{1}{\Gamma(X_1,p)} + \frac{1}{\Gamma(X_2,p)} + \alpha\right)^{-1} \le \Gamma(X_1cX_2,p) \le \left(\frac{1}{\Gamma(X_1,p)} + \frac{1}{\Gamma(X_2,p)}\right)^{-1},
\end{equation*}
with
\begin{equation*}
\alpha = \frac{1}{2}\cdot\frac{(\lVert u_p^1\rVert_1 + \lVert u_p^2\rVert_1)^2}{2c^p - M_p(X_1) - M_p(X_2)}
\end{equation*}

Note that $M_p(\{x\}) = 0$ and $\Gamma(\{x\},p) = \infty$, so $\Gamma(\{x\},p)^{-1} := 0$.
}
\proof{It is shown in Theorem 3.5 in \cite{RW} that $\Gamma(X,p) = \frac{2}{\beta(X,p)}$, with $\beta(X,p) = \\
\underset{z \in\{-1,1\}^n}{\max}\widehat{(D_p(X)}z|z)$, for $(X,d)$ a finite metric space of strict $p$-negative type of at least two points. Applying Corollary \ref{Cor66} we are done.
\begin{flushright}$\qedsymbol$\end{flushright}
}

\section{An application to finite ultrametric spaces}\label{7}

Let $(X,d)$ be a metric space. Recall that $d$ is an ultrametric on $X$, ($(X,d)$ is called an ultrametric space) if $d(x,y) \le \max(d(x,z),d(z,y))$ for all $x, y, z$ in $X$. Further recall that a real symmetric $n \times n$ matrix $A = (a_{ij})_{i,j = 1}^n$ with nonnegative entries is called strictly ultrametric if

\begin{enumerate}
\item $a_{ij} \ge \min\{a_{ik},a_{kj}\}$, for all $i, j, k$ and 
\item $a_{ii} > \max\{a_{ij}, j \neq i\}$, for all $i$
\end{enumerate}
holds.

The following result was shown in \cite{TF}:

\lem\label{Lem71}{Let $(X,d)$ be a finite ultrametric space. Then we have
\begin{enumerate}
\item $(X,d^p)$ is a finite ultrametric space, for each $p > 0$
\item $(X,d)$ is of strict $p$-negative type, for each $p > 0$.
\end{enumerate}
}

\lem\label{Lem72}{Let $(X,d)$, $|X| = n$ be a finite ultrametric space with $p$-distance matrix $D_p(X)$ and diameter $\Delta(X)$. Then we have $\Delta(X)^p\1\1^T - D_p(X)$ is strictly ultrametric for each $p > 0$.
}
\proof{Follows by definition and Lemma \ref{Lem71} part 1.}

\lem\label{Lem73}{Let $(X,d)$ be a finite ultrametric of at least two points. Then we have
\begin{enumerate}
\item $D_p(X)$ is nonsingular
\item The entries of $D_p(X)^{-1}\1$ are all positive.
\end{enumerate}
}
\proof{By Lemma \ref{Lem71} part 2, we know that $D_p(X)^{-1}$ exists. It is shown in \cite{SM} (Lemma 1) that for each strictly ultrametric $n \times n$ matrix $A$ there exists some $x_0$ in $\R^n$, such that $Ax_0 = \1$ and $(x_0|e_i) > 0$, for all $1 \le i \le n$.\\
Applying Lemma \ref{Lem72} we get some $x_0$ in $\R^n$ with $(\Delta(X)^p\1\1^T - D_p(X))x_0 = \1$. Hence $D_p(X)x_0 = (\Delta(X)^p(x_0|\1) - 1)\1$. Since $(x_0|e_i) > 0$ for all $1 \le i \le n$ we get $(D_p(X)x_0|e_i) > 0$ for all $1 \le i \le n$ and hence $\Delta(X)^p(x_0|\1) -1 > 0$. Therefore all entries of $D_p(X)^{-1}\1 = (\Delta(X)^p(x_0|\1) - 1)^{-1}x_0$ are positive.
\begin{flushright}$\qedsymbol$\end{flushright}
}

\lem\label{Lem74}{Let $(X,d)$, $|X| = n$ be a finite ultrametric space. Then we have\\
$M_p(X) \le \frac{n - 1}{n}\Delta(X)^p$, for each $p > 0$.}
\proof{Since $M_p(\{x\}) = 0$ we may assume that $n \ge 2$. Let $u_p$ be the unique element in $F_1$, such that $D_p(X)u_p = M_p(X)\1$ (see Remark \ref{Rem310}).
Now $D_p(X)^{-1}\1 = M_p(X)^{-1}u_p$. Let $u_p^T = (\alpha_1, \alpha_2, \ldots, \alpha_n)$. By Lemma \ref{Lem73} part 2 we have $\alpha_i > 0$, for all $1 \le i \le n$ and $\alpha_1 + \alpha_2 + \ldots + \alpha_n = 1$, since $u_p$ in $F_1$. Therefore

\begin{equation*}
\begin{multlined}
M_p(X) = (D_p(X)u_p|u_p) = \sum_{i,j = 1}^nd(x_i,x_j)^p\alpha_i\alpha_j \le\\
\le \Delta(X)^p\sum_{i = 1}^n\alpha_i(1 - \alpha_i) = \Delta(X)^p(1 - \sum_{i = 1}^n\alpha_i^2).
\end{multlined}
\end{equation*}

Since
\begin{equation*}
\sum_{i = 1}^n\alpha_i^2 \ge \frac{1}{n}\left(\sum_{i = 1}^n\alpha_i\right)^2 = \frac{1}{n},
\end{equation*}
we are done.
\begin{flushright}$\qedsymbol$\end{flushright}

The following assertion is clear by definition of ultrametric spaces:
\rem\label{Rem75}{Let $(X_1,d_1)$ and $(X_2,d_2)$ be two finite ultrametric spaces, such that $X_1 \cap X_2 = \emptyset$ and let $c$ be a real number with $c \ge \max(\Delta(X_1),\Delta(X_2))$. Then we have: $X_1cX_2$ is a finite ultrametric space.}

The next result is a kind of converse to Remark \ref{Rem75} and was shown in order to decompose strictly ultrametric metrices (see Propostion 2.1 in \cite{RN}).

We translate the proof given in \cite{RN} into the language of ultrametric spaces:

\thm\label{Thm76}{Let $(X,d)$ be a finite ultrametric space of at least two points. Then there are subspaces $X_1, X_2$ of $X$ such that
\begin{enumerate}
\item $X_1 \neq \emptyset$, $X_2 \neq \emptyset$, $X_1 \cap X_2 = \emptyset$, $X_1 \cup X_2 = X$
\item $X = X_1\Delta(X)X_2$
\end{enumerate}
}

\proof{Take two points $x_1 \neq x_2$ in $X$ such that $d(x_1,x_2) = \Delta(X)$. Let\\
$X_1 = \{x \in X | d(x,x_1) < \Delta(X)\}$ and $X_2 = \{x \in X | d(x,x_1) = \Delta(X)\}$.
Since $x_1 \in X_1$ and $x_2 \in X_2$ we get $X_1 \neq \emptyset$ and $X_2 \neq \emptyset$. Of course we have $X_1 \cap X_2 = \emptyset$ and $X_1 \cup X_2 = X$.
Now let $z_1, z_2$ be two points in $X$ with $z_1 \in X_1$ and $z_2 \in X_2$.

$\Delta(X) = d(z_2,x_1) \le \max(d(z_2,z_1),d(z_1,x_1))$. Since $d(z_1,x_1) < \Delta(X)$ we get $d(z_1,z_2) = \Delta(X)$. Hence $X = X_1\Delta(X)X_2$.
\begin{flushright}$\qedsymbol$\end{flushright}
}

\thm\label{Thm77}{Let $(X,d)$ be a finite ultrametric space of at least two points. Let $X = X_1\Delta(X)X_2$ acording to Theorem \ref{Thm76}. Then we have
\begin{equation*}
\left(\frac{1}{\Gamma(X_1,p)} + \frac{1}{\Gamma(X_2,p)} + \alpha\right)^{-1} \le \Gamma(X,p) \le \left(\frac{1}{\Gamma(X_1,p)} + \frac{1}{\Gamma(X_2,p)}\right)^{-1}
\end{equation*}
with
\begin{equation*}
\alpha = \frac{2}{2\Delta(X)^p - \frac{|X_1| - 1}{|X_1|}\Delta(X_1)^p - \frac{|X_2| - 1}{|X_1|}\Delta(X_2)^p} \le \frac{|X|}{\Delta(X)^p},
\end{equation*}
where
$\Gamma(\{x\},p)^{-1} = 0$.
}

\proof{Applying Theorem \ref{Thm67} it remains to show that
\begin{equation*}
\frac{1}{2}\cdot\frac{(\lVert u_p^1\rVert_1 + \lVert u_p^2\rVert_1)^2}{2\Delta(X)^p - M_p(X_1) - M_p(X_2)} \le \alpha.
\end{equation*}

For one point spaces we have $u_p = 1$ and $M_p(\{x\}) = 0$. If $X_i$ has at least two points, we apply Lemma \ref{Lem73} and obtain that all entries of $u_p^i = M_p(X_i)D_p(X_i)^{-1}\1$ are positive and hence $\lVert u_p^{i}\rVert_1 = 1$, since $u_p^{i}$ in $F_1$. Lemma \ref{Lem74} implies $M_p(X_i) \le \frac{|X_i| - 1}{|X_i|}\Delta(X_i)^p$.

Summing up we get $\lVert u_p^1\rVert_1 = \lVert u_p^2\rVert_1 = 1$ and $2\Delta(X)^p - M_p(X_1) - M_p(X_2) \ge 2\Delta(X)^p - \frac{|X_1| - 1}{|X_1|}\Delta(X_1)^p - \frac{|X_2| - 1}{|X_2|}\Delta(X_2)^p$.

Hence
\begin{equation*}
\begin{multlined}
\frac{1}{2}\cdot\frac{(\lVert u_p^1\rVert_1 + \lVert u_p^2\rVert_1)^2}{2\Delta(X)^p - M_p(X_1) - M_p(X_2)} \le \\
\le \frac{2}{2\Delta(X)^p - \frac{|X_1| -1}{|X_1|}\Delta(X_1)^p - \frac{|X_2| - 1}{|X_2|}\Delta(X_2)^p} = \alpha \le \frac{|X|}{\Delta(X)^p}.
\end{multlined}
\end{equation*}
\begin{flushright}$\qedsymbol$\end{flushright}
}

We now estimate $\Gamma(X,p)$ for a concrete finite ultrametric space acording to Theorem \ref{Thm77}. First recall that for a given finite connected edge-weighted graph $G = (V,E,w)$, $w:E \rightarrow \R^+$, $d : V \times V \rightarrow \R$, $d(u,v) = \min\{\bar{w}(W), W$ a walk connecting $u$ and $v\}$ with $\bar{w}(W) = \max\{w(e), e$ a edge on $W\}$ the space $(V,d)$ is a finite ultrametric space.

\exa{Let $G = (V,E,w)$ be:

$V = \{a,b,c,d,e,f,g\},\\
E = \{\{a,b\},\{b,c\},\{c,d\},\{c,e\},\{e,f\},\{f,g\}\},\\
w(\{c,d\}) = w(\{e,f\}) = 1\\
w(\{a,b\}) = w(\{b,c\}) = 2\\
w(\{c,e\}) = 3\\
w(\{f,g\}) = 4$.

The distance matrix of $(V,d)$ is 

\begin{equation*}
\begin{pmatrix}
0 & 2 & 2 & 2 & 3 & 3 & 4 \\ 2 & 0 & 2 & 2 & 3 & 3 & 4 \\ 2 & 2 & 0 & 1 & 3 & 3 & 4 \\ 2 & 2 & 1 & 0 & 3 & 3 & 4 \\ 3 & 3 & 3 & 3 & 0 & 1 & 4 \\ 3 & 3 & 3 & 3 & 1 & 0 & 4 \\ 4 & 4 & 4 & 4 & 4 & 4 & 0
\end{pmatrix}.
\end{equation*}

Now we get $V = \{a,b,c,d,e,f,g\} = \{a,b,c,d,e,f\} 4 \{g\}$\\
$\{a,b,c,d,e,f\} = \{a,b,c,d\} 3 \{e,f\}$\\
$\{a,b,c,d\} = \{a,b\} 2 \{c,d\}$ by Theorem \ref{Thm76}.

Let $V_1 = \{a,b,c,d,e,f\}$, $V_2 = \{g\}$, $V_3 = \{a,b,c,d\}$, $V_4 = \{e,f\}$, $V_5 = \{a,b\}$, $V_6 = \{c,d\}$.

Applying Theorem \ref{Thm77} we obtain

\begin{equation*}
\begin{multlined}
\frac{1}{\Gamma(V_1,p)} \le \frac{1}{\Gamma(V,p)} \le \frac{1}{\Gamma(V_1,p)} + \frac{7}{4^p} \\
\frac{1}{\Gamma(V_3,p)} + \frac{1}{\Gamma(V_4,p)} \le \frac{1}{\Gamma(V_1,p)} \le \frac{1}{\Gamma(V_3,p)} + \frac{1}{\Gamma(V_4,p)} + \frac{6}{3^p} \\
\frac{1}{\Gamma(V_5,p)} + \frac{1}{\Gamma(V_6,p)} \le \frac{1}{\Gamma(V_3,p)} \le \frac{1}{\Gamma(V_5,p)} + \frac{1}{\Gamma(V_6,p)} + \frac{4}{2^p}
\end{multlined}
\end{equation*}
Summing up we get

\begin{equation*}
\begin{multlined}
\frac{1}{\Gamma(V_4,p)} + \frac{1}{\Gamma(V_5,p)} + \frac{1}{\Gamma(V_6,p)} \le \frac{1}{\Gamma(V,p)} \le \\
\frac{1}{\Gamma(V_4,p)} + \frac{1}{\Gamma(V_5,p)} + \frac{1}{\Gamma(V_6,p)} + \frac{4}{2^p} + \frac{6}{3^p} + \frac{7}{4^p}
\end{multlined}
\end{equation*}

Since $\Gamma(V_5,p) = \Gamma(2\cdot X_2,p) = 2^p\Gamma(X_2,p) = 2^p$, $X_2$ the discrete space consisting of two points, we get
\begin{equation*}
\underset{p \rightarrow \infty}{\lim}\frac{1}{\Gamma(V,p)} = \frac{1}{\Gamma(V_4,p)} + \frac{1}{\Gamma(V_6,p)},
\end{equation*}
where $V_4$ and $V_6$ are the maximal (with respect to inclusion) discrete subspaces of at least two points of $V$.}

The last limit equation is not a coincidence:

Recently the following beautiful formula was shown by Doust, S\'{a}nchez and Weston.

\thm\label{Thm79}{(=Theorem 1.5 in \cite{ID3}) Let $(X,d)$ be a finite ultrametric space with at least two points and minimum non-zero distance $\alpha$. If $B_1, B_2, \ldots, B_e$ are the distinct coteries of $(X,d)$ then
\begin{equation*}
\underset{p \rightarrow \infty}{\lim}\frac{\Gamma(X,p)}{\alpha^p} = \left(\frac{1}{\gamma(|B_1|)} + \ldots + \frac{1}{\gamma(|B_e|)}\right)^{-1}
\end{equation*}
with
\begin{equation*}
\gamma(n) = \frac{1}{2}\left(\Big\lfloor \frac{n}{2}\Big\rfloor^{-1} + \Big\lceil \frac{n}{2}\Big\rceil^{-1}\right)
\end{equation*}
and a coterie of $(X,d)$ is defined as $\{x \in X | d(x,y) \le \alpha\}$, for some $z$ in $X$ with\\
$|\{x \in X | d(x,y) \le \alpha\} | \ge 2$.

\rem{Let $(X,d)$ be a finite ultrametric space of at least two points and the distinct distances in $(X,d)$ are given by $0 = \alpha_0 < \alpha_1 < \alpha_2 < \ldots < \alpha_k = \Delta(X)$. We assume $1 = \alpha_1 < \Delta(X)$ \quad ($\frac{\Gamma(X,p)}{\alpha_1^p} = \Gamma\left(\frac{1}{\alpha_1}\cdot X,p\right)$).

Applying Theorem \ref{Thm77} a several times, we get
\begin{equation*}
0 \le \frac{1}{\Gamma(X,p)} - \left(\sum_{i = 1}^{e}\frac{1}{\gamma(|B_i|)}\right) \le \frac{A}{\alpha_2^p},
\end{equation*}
where $A > 0$ is independent of $p$ and $B_1, \ldots, B_e, \gamma$ defined as in Theorem \ref{Thm79} ($B_i = X_{|B_i|}$ and $\gamma(|B_i|) = \Gamma(X_{|B_i|},p)$, for all $p > 0$).

Hence
\begin{equation*}
\underset{p \rightarrow \infty}{\lim}\frac{1}{\Gamma(X,p)} = \sum_{i = 1}^{e}\frac{1}{\gamma(|B_i|)}
\end{equation*}
and the rate of convergence is at least $\alpha_2^{-p}$.

\end{document}